\documentclass[12pt,a4paper,reqno]{amsart}
\usepackage{amssymb}
\usepackage{dcpic,pictex}
\usepackage{multirow}
\usepackage{graphicx}
\usepackage{slashbox}
\usepackage{cite}

\usepackage{hyperref}
\newtheorem{Def}{Definition}[section]

\newtheorem{The}{Theorem}[section]
\newtheorem{Rem}{Remark}[section]
\newcommand{\R}{\mathbb{R}}

\newcommand{\N}{\mathbb{N}}

\begin{document}

%\selectlanguage{english}

\title[Fractional derivatives and the fundamental theorem]{Fractional derivatives and the fundamental theorem of Fractional Calculus}

\author{Yuri Luchko}
%\address{Beuth University of Applied Sciences Berlin}
\curraddr{Beuth Technical University of Applied Sciences Berlin,  
     Department of  Mathematics, Physics, and Chemistry,  
     Luxemburger Str. 10,  
     13353 Berlin,
Germany}
\email{luchko@beuth-hochschule.de}

%    General info
%\UDC{519.21}
%\date{04.07.2020}
\dedicatory{}
% \thanks{Paper accepted for publication in {\it Fract. Calc. Appl. Anal.}, to appear at https://www.degruyter.com/view/j/fca}
% published in 
% Fract. Calc. Appl. Anal.,Vol. 23, No 4 (2020), pp. 939–966, DOI: 10.1515/fca-2020-0049

\begin{abstract}
In this paper, we address the one-parameter families of the fractional integrals and derivatives defined on a finite interval. First we remind the reader of the known fact that under some reasonable conditions, there exists precisely one unique family of the fractional integrals, namely, the well-known Riemann-Liouville fractional integrals. As to the fractional derivatives, their natural definition follows from the fundamental theorem of the Fractional Calculus, i.e., they are introduced as the left-inverse operators to the Riemann-Liouville fractional integrals. Until now, three families of such derivatives were suggested in the literature: the Riemann-Liouville fractional derivatives, the Caputo fractional derivatives, and the Hilfer fractional derivatives. We clarify the interconnections between these derivatives on different spaces of functions and provide some of their properties including the formulas for their projectors and the Laplace transforms. However, it turns out that there exist infinitely many other families of the fractional derivatives that are the left-inverse operators to the Riemann-Liouville fractional integrals. In this paper, we focus on an important class of these fractional derivatives and discuss some of their properties. 
\end{abstract}

\maketitle

%%%%%%%%%%%%%%%%%%%%%%%%%%%%%%%%%%%%%%%%%%%%%%%%
\section{Introduction} \label{sec:1}

\setcounter{section}{1}
\setcounter{The}{0}

Within the last few years, a lot of efforts of the Fractional Calculus (FC) community was put into clarifying the question what are the fractional integrals and derivatives and what are they not (\cite{DGGS, Giu18, Han20, Hil19, HilLuc, Ort15, Ort18, Sty18, Tar13, Tar19}). These discussions mainly concerned the ``new fractional integrals and derivatives" (\cite{DGGS, Giu18, Han20, Ort18, Sty18, Tar13, Tar19}), the integro-differential operators of convolution type with some general kernels (\cite{DGGS, Han20, Koch11, Koch19, Sty18, Zac08, Zac09}), and the abstract axioms of FC
(\cite{HilLuc, Ort15, Ross75}).

In this framework, the ``classical" definitions of the fractional integrals and derivatives as the Riemann-Liouville integral and derivative (\cite{KocLuc19, Samko}), the Caputo derivative (\cite{Cap, Diet10, KocLuc19}), the generalized Riemann-Liouville  derivative or the Hilfer derivative (\cite{Hil00, Hil, KocLuc19}), etc., are usually considered to be postulated. In this paper, we are going to take a closer and critical look at the ``right" definitions of the one-parameter families of the fractional integrals and derivatives defined on a finite interval.

In the 2nd section, we reproduce a result derived in \cite{Cart78} more than forty years ago  regarding the ``right" one-parameter families of the fractional integrals  defined on a finite interval. It turns out that under some reasonable conditions, the only family of the fractional integrals on a finite interval are the Riemann-Liouville fractional integrals.

The question regarding the ``right" fractional derivatives is more delicate and is addressed in the 3rd section. In calculus, the integer order derivatives are usually defined via the limits of the difference quotients. The definite integral is introduced independently through the Riemann sums or in the Lebesgue sense. The fundamental theorem of calculus establishes a connection between these two independently defined objects and says that the first order derivative is a left-inverse operator to the definite integral with the variable upper limit of integration on a suitable space of functions. Should we follow this approach to introduce the ``right" fractional derivative, we naturally arrive at the Gr\"unwald-Letnikov definition (\cite{Diet10, Grue, Let, Samko}). This definition is useful and important, say,  for  numerical calculation of the fractional derivatives. However, it is  unpractical for analytical investigations and usually replaced by other definitions that are equivalent to the Gr\"unwald-Letnikov definition on some suitable spaces of functions. Say, for the functions from $C^n[a,b]$ with $n-1<\alpha \le n,\ n \in \N$, the 
Gr\"unwald-Letnikov fractional derivative of order $\alpha$ coincides with the Riemann-Liouville fractional derivative of order $\alpha$ (\cite{Diet10}). That's why in analysis of the operator-theoretic properties of the fractional derivatives, in the fractional differential equations, etc., the Gr\"unwald-Letnikov derivative is often replaced by the Riemann-Liouville  derivative.

Another - and in fact the standard - approach for introducing the fractional derivatives is to assume that they are connected to the corresponding fractional integrals by a fractional analogy of the fundamental theorem of calculus, i.e., to define them as the left-inverse operators to the fractional integrals. As already said, the only ``right"  fractional integrals on a finite interval are the Riemann-Liouville integrals $I^\alpha,\ \alpha\ge 0$ and thus in the framework of this approach, the Abel integral equation plays a decisive role. It was known already to Abel (\cite{Abel1, Abel2}), that - under some natural conditions - the Abel integral equation has a unique solution and this solution is given by the formula that is presently known as the Riemann-Liouville fractional derivative (see \cite{Samko} for a detailed derivation of the solution formula). Thus, on the function space  $I^\alpha(L_1(a,b))$, there is only one ``right" one-parameter family of the fractional derivatives, namely, the Riemann-Liouville fractional derivatives.

One of the most important and powerful ideas in mathematics in general and in FC in particular  is extension of validity domains for formulas, statements, functions, operators, etc. Following this strategy, the basic definition domain $I^\alpha(L_1(a,b))$ of the Riemann-Liouville fractional derivative (as solution formula to the Abel integral equation) is extended to larger spaces of functions. Say, in the case $0<\alpha < 1$, to the space of functions, whose Riemann-Liouville integrals of order $1-\alpha$ are absolutely continuous functions on the interval $[a,\, b]$. This extended operator is also called the Riemann-Liouville fractional derivative even if - strictly speaking - any operator is uniquely defined not only by its values, but also by its definition domain. The extended Riemann-Liouville fractional derivative of order $\alpha$ keeps the property to be a left-inverse operator to the Riemann-Liouville fractional integral on even larger spaces of functions, e.g., on $L_1(a,b)$ (\cite{Samko}).

At the first look, the considerations described  above let no place for other fractional derivatives on a finite interval, say, for the Caputo or the Hilfer derivatives. However, it is known that these derivatives are also the left-inverse operators to the Riemann-Liouville fractional integrals on some suitable spaces of functions  (\cite{Diet10, Hil}) and thus they belong to the class of the fractional derivatives in the sense of the fundamental theorem of FC. As we will see in Section \ref{sec:3}, % the 3rd section, %
the construction of these derivatives (especially of their definition domains) is not as straightforward as the one of the Riemann-Liouville fractional derivative. Namely, instead of extension one has to start with a contraction of the basic space  $I^\alpha(L_1(a,b))$ (very unusual in mathematics!) and to consider its subspaces, where the first order derivative commutes with the corresponding Riemann-Liouville fractional integrals.  Only after this, in the second step, the definition domains of the Caputo or the Hilfer fractional derivatives are extended from these subspaces to the larger spaces of functions.  It is important to note that on their extended definition domains these derivatives do not coincide anymore with the Riemann-Liouville fractional derivative. In particular, their kernels are different and this means that one has to provide different initial conditions while considering the initial-value problems for the fractional differential equations with these fractional derivatives.

The construction of the Hilfer fractional derivative can be extended to the finite chains of compositions of the first order derivatives and some suitable Riemann-Liouville fractional integrals. We call these operators the
$n$-th level fractional derivatives (the Riemann-Liouville, the Caputo, and the Hilfer derivatives can be interpreted as the 1st level fractional derivatives). It turns out that the $n$-th level fractional derivatives are the left-inverse operators to the Riemann-Liouville fractional integrals on some suitable spaces of functions and thus satisfy the fundamental theorem of FC. As a result, there exist infinitely many different ``right" one-parameter families of the fractional derivatives defined on a finite interval. All of them are the left-inverse operators to the Riemann-Liouville fractional integrals on some suitable spaces of functions.

Finally, in Section \ref{sec:4}, we discuss some important properties of the fractional derivatives mentioned above. In particular, for the 2nd level fractional derivative, we derive an explicit formula for its projector that determines the number and the form of the initial conditions for the Cauchy problems for the fractional differential equations with this fractional derivative. It is worth mentioning that the kernel dimension of the $n$-th level fractional derivative of order $\alpha,\ 0<\alpha <1$ can be equal to or less than $n$ and thus, in the general case, up to $n$ initial conditions are required to guarantee uniqueness of solutions  to the Cauchy problems for the fractional differential equations with this derivative. Then we shortly discuss definition and properties  of the 2nd level  fractional derivatives on the positive real semi-axis and derive a formula for its Laplace transform.

\vspace*{-3pt}
%%%%%%%%% Section 2 %%%%%%%%%%%%%%%%%%%
\section{Fractional integrals on a finite interval} \label{sec:2}

\setcounter{section}{2}
\setcounter{The}{0}

We start this section with formulation of a conjecture that was posed by J.S. Lew during the first FC conference  held in New Haven, USA in 1974 (\cite{Osl}).

Let $L$  be either $L_1(0,1)$ or $L_2(0,1)$. The  formula for the Riemann-Liouville fractional integrals ($x\in[0,\, 1]$)
\vskip -10pt
\begin{equation}
\label{RLI}
(I^\alpha\, f)(x) = \begin{cases}
\frac{1}{\Gamma(\alpha)} \int_0^x (x-t)^{\alpha -1}\, f(t)\, dt, & \alpha >0,\\
f(x), & \alpha = 0,
\end{cases}
\end{equation}
yields a family $I^\alpha,\ \alpha \ge 0$ of the bounded linear operators on $L$ satisfying the properties:

\vspace{0.2cm}

\noindent
L1) $(I^0 \, f)(x) = f(x),\ f\in L$ (1st interpolation condition),
\vspace{0.1cm}

\noindent
L2) $(I^1\, f)(x) = \int_0^x f(t)\, dt, \ f\in L$ (2nd interpolation condition),
\vspace{0.1cm}

\noindent
L3) $(I^\alpha \, I^\beta\, f)(x) = (I^{\alpha + \beta}\, f)(x), \ \alpha,\beta \ge 0, \ f\in L$ (index law),
\vspace{0.1cm}

\noindent
L4) $I^\alpha$ is continuous in $\alpha$ in  some operator topology (continuity),
\vspace{0.1cm}

\noindent
L5) $f\in L$ and $f(x)\ge 0$ a.e. $\Rightarrow$ $(I^\alpha\, f)(x) \ge 0$ a.e. for all $\alpha \ge 0$ (non-negativity).

\vspace{0.2cm}

Lew conjectured that the properties L1)-L5) uniquely determine the family of the Riemann-Liouville fractional integrals defined by \eqref{RLI}.

Two years later, in 1976, the conjecture by Lew
was confirmed  by Cartwright and  McMullen in \cite{Cart78} (\cite{Cart78} was published in 1978, but submitted in 1976). Moreover, they proved the uniqueness of the family of the Riemann-Liouville fractional integrals under weaker conditions than those formulated by Lew and for larger spaces of functions. Their main result is as follows:

\vspace*{-4pt}

\begin{The}[\cite{Cart78}]
\label{t_1}
Let $E$ be the space $L_p(0,1),\, 1\le p <+\infty$, or $C[0,1]$. Then there is precisely one family $I_\alpha,\, \alpha >0$ of operators on $E$ satisfying the following conditions:

\vspace{0.2cm}

\noindent
CM1)  $(I_1\, f)(x) = \int_0^x f(t)\, dt, \ f\in E$ (interpolation condition),
\vspace{0.1cm}

\noindent
CM2) $(I_\alpha \, I_\beta\, f)(x) = (I_{\alpha + \beta}\, f)(x), \ \alpha,\beta > 0, \  f\in E$ (index law),
\vspace{0.1cm}

\noindent
CM3) $\alpha \to I_\alpha$ is a continuous map of $(0,\, +\infty)$ into ${\mathcal L}(E)$ for some Hausdorff topology on ${\mathcal  L}(E)$, weaker than the norm topology (continuity),

\noindent
CM4) $f\in E$ and $f(x)\ge 0$ (a.e. for $E=L_p(0,1)$) $\Rightarrow$ $(I_\alpha\, f)(x) \ge 0$ (a.e. for $E=L_p(0,1)$) for all $\alpha > 0$ (non-negativity).

\vspace{0.2cm}

\noindent
That family is given by the Riemann-Liouville formula \eqref{RLI} with $\alpha >0$.
\end{The}

Please note that the result by Cartwright and  McMullen does not involve the case $\alpha =0$ and the 1st interpolation condition L1) from the conjecture by Lew. However, the family of the Riemann-Liouville fractional integrals defined for $\alpha >0$ can be uniquely extended to the family defined for $\alpha \ge 0$ by setting $I^0$ to be the identity operator. This  is justified by the following statement (Theorem 2.6 in \cite{Samko}): The family of the Riemann-Liouville fractional integrals $I^\alpha,\, \alpha \ge 0$ given by \eqref{RLI} forms a semigroup in $L_p(0,\ 1),\ p\ge 1$, which is strongly continuous for all $\alpha \ge 0$, i.e., the relation
 \vskip -13pt
\begin{equation}
\label{cont}
\lim_{\alpha \to \alpha_0} \| I^\alpha f - I^{\alpha_0} f\|_{L_p(0,1)} = 0
\end{equation}
\vskip -3pt \noindent
holds valid for any $\alpha_0,\ 0\le \alpha_0 < +\infty$ and for any $f\in L_p(0,1)$. Evidently, this extended family of operators satisfies the properties L1)-L5) for all $\alpha \ge 0$.

%%%%%%

From the present viewpoint (\cite{HilLuc}), the conditions CM1)-CM4) are very natural and in fact desirable for any definition of the fractional integrals defined on a finite interval. As proved in \cite{Cart78}, they are also sufficient for uniqueness of the family of the Riemann-Liouville fractional integrals. Thus, in this sense, the Riemann-Liouville fractional integrals are the only ``right" one-parameter fractional integrals defined on a finite interval.

In the further discussions, we need the well-known formula for the Riemann-Liouville fractional integral of a power law function:
\begin{equation}
\label{power}
(I^\alpha\, t^{\beta})(x) = \frac{\Gamma(\beta +1)}{\Gamma(\alpha+\beta +1)}\, x^{\alpha + \beta},\ \
\alpha \ge 0,\ \beta > -1.
\end{equation}
For a very detailed discussion of other properties of the Riemann-Liouville fractional integrals, we refer to \cite{Samko}.

\vskip 2pt

In what follows, without loss of generality, we restrict ourselves to the operators defined on the interval $[0,\, 1]$ (the case of the interval $[a,\ b]$ can be reduced to the case of the interval $[0,\, 1]$ by a linear variables substitution).

In the rest of this section, we recall some  basic facts regarding the Abel integral equation
\vskip -13pt
\begin{equation}
\label{Abel}
(I^\alpha\, \phi)(x) = \frac{1}{\Gamma(\alpha)} \int_0^x (x-t)^{\alpha -1}\, \phi(t)\, dt\ = f(x),\ \,
0<\alpha<1,\ x\in[0,\, 1].
\end{equation}
In \cite{Abel1, Abel2}, Abel solved this equation by applying the operator $I^{1-\alpha}$ to both sides of the equation \eqref{Abel} and by using the semigroup property of the Riemann-Liouville fractional integrals (in modern terminology):
\vskip -10pt
$$
(I^{1-\alpha}\, I^\alpha\, \phi)(x) = (I^1\, \phi)(x) = \int_0^x \phi(t)\, dt = (I^{1-\alpha}\, f)(x).
$$
To determine the unknown function $\phi$, he differentiated the last formula and applied the fundamental theorem of calculus ($x\in[0,\, 1])$:
\vskip -10pt
\begin{equation}
\label{Asol}
\phi(x)  = \frac{d}{dx}\, (I^{1-\alpha}\, f)(x) = \frac{d}{dx}\, \frac{1}{\Gamma(1-\alpha)} \int_0^x (x-t)^{-\alpha}\, f(t)\, dt,\ \, 0<\alpha<1.
\end{equation}

In the next section, we will employ the following result regarding solvability of the Abel integral equation in $L_1(0,1)$:

\vspace*{-3pt}

\begin{The}[\cite{Samko}]
\label{t_A}
The Abel integral equation \eqref{Abel} is solvable in $L_1(0,1)$ if and only if
\begin{equation}
\label{Acond}
I^{1-\alpha}\, f\in \mbox{AC}([0,1]) \ \ \mbox{and} \ \ (I^{1-\alpha}\, f)(0) = 0.
\end{equation}
If these conditions are satisfied, the Abel integral equation has a unique solution in $L_1(0,1)$ given by the formula \eqref{Asol}.
\end{The}

The notation $\mbox{AC}([0,1])$ stands for the space of functions that are absolutely continuous on the interval $[0,\, 1]$. This space can be characterized as follows:
\vskip -10pt
\begin{equation}
\label{absc}
f\in \mbox{AC}([0,1]) \, \Leftrightarrow \, \exists \phi \in L_1(0,1):\, f(x) = f(0) + \int_0^x \phi(t)\, dt,\ x\in [0,\, 1].
\end{equation}

\vspace*{-2pt}

\begin{Rem}
\label{r_0}
Based on the representation \eqref{absc}, a (weak) derivative of a function $f\in \mbox{AC}([0,1])$ can be defined:
\vskip -10pt
\begin{equation}
\label{absc_d}
f(x) = f(0) + \int_0^x \phi(t)\, dt,\ x\in [0,\, 1] \ \Rightarrow \ \frac{df}{dx} := \phi \in L_1(0,1).
\end{equation}
In the further discussions, we interpret the first order derivative of the absolutely continuous functions in the sense of the formula  \eqref{absc_d}.
\end{Rem}

\vspace*{-4pt}
%%%%%%%%%%%% Section 3 %%%%%%%%%%%%%%%%%%%%%%%%%%%%%%%%%%%%
\section{Fractional derivatives on a finite interval} \label{sec:3}

\setcounter{section}{3}
\setcounter{The}{0}

\vskip 2pt 

As discussed in Introduction, the most standard and natural way for defining the fractional derivatives is by means of their connection to the fractional integrals. It turns out that the definitions and the properties of the fractional derivatives essentially depend on the spaces of functions, where they are defined (an operator is a triple $(A,X,Y)$ consisting of the domain $X$, the range $Y$, and the correspondence $A:\ X\to Y$).  In this section, we discuss both the known fractional derivatives like the Riemann-Liouville, the Caputo, and the Hilfer fractional derivatives and the new fractional derivatives that we call the $n$-th level fractional derivatives. First we introduce a general notion of a one-parameter family of the fractional derivatives defined on a finite interval.

\vspace*{-3pt}

\begin{Def}
\label{def1}
Let $I^\alpha,\ \alpha \ge 0$ be the family of the fractional Riemann-Liouville integrals defined by \eqref{RLI}. A one-parameter family $D^\alpha$, $\alpha \ge 0$ of the linear operators is called the fractional derivatives if and only if it satisfies the Fundamental Theorem of Fractional Calculus formulated below.
\end{Def}

\vspace*{-4pt}

\begin{The}[Fundamental Theorem of FC]
\label{t2}
For the fractional derivatives $D^\alpha,\ \alpha \ge 0$ and the Riemann-Liouville fractional  integrals $I^\alpha,\ \alpha \ge 0$, the relation
\vskip -12pt
\begin{equation}
\label{ftFC}
(D^\alpha\, I^\alpha\, \phi)(x) = \phi(x),\ x\in [0,\, 1]
\end{equation}
\vskip 4pt \noindent
holds true on appropriate nontrivial spaces of functions.
\end{The}

In fact, the Fundamental Theorem of FC is a part of Definition \ref{def1}, i.e., a linear operator is called a fractional derivative if and only if it is a left-inverse operator to the Riemann-Liouville fractional integral on a certain space of functions (compare to the desiderata (c) from \cite{HilLuc}). It turns out that there exist infinitely many different families of the fractional derivatives in the sense of Definition \ref{def1}. In the rest of this section, we discuss some known and new families of the fractional derivatives on the interval $[0,\, 1]$.

\vspace*{-3pt}

\begin{Rem}
\label{r_ft}
In calculus, the formula of type \eqref{ftFC} with $\alpha = 1$ is usually called the 1st fundamental theorem of calculus. The 2nd fundamental theorem of calculus states that $\int_{0}^x\, f^\prime(t)\, dt \, = \, f(x)-f(0)$. We address the FC analogy of the 2nd fundamental theorem in Section \ref{sec:4}.
\end{Rem}

\vspace*{-15pt}

\begin{Rem}
\label{r_1}
It is worth mentioning that the formula \eqref{ftFC} and the relation $(I^0\, f)(x) = f(x)$ uniquely define the fractional derivative of order $\alpha =0$ as the identity operator: $(D^0\, f)(x) = f(x)$. Therefore, in what follows, we mainly restrict our attention to the case $\alpha >0$.
\end{Rem}

\vspace*{-15pt}

\begin{Rem}
\label{r_2}
The Riemann-Liouville fractional integral  $I^\alpha$ is injective on $L_1(0,1)$, i.e., its kernel contains only the null-function $f(x)=0$ a.e. on $[0,\, 1]$ (Lemma 2.5 in \cite{Samko}). Evidently, this statement holds true for any linear operator that possesses a linear left-inverse operator. For $I^\alpha$, it follows from the realization \eqref{ftFC_RL} of the Fundamental Theorem of FC for the Riemann-Liouville fractional derivatives.
\end{Rem}

\vspace*{-15pt}

\begin{Rem}
\label{r_3}
For the sake of formulations simplicity, in what follows, we restrict ourselves to the function space  $L_1(0,1)$ and its subspaces (a similar theory can be developed for, say, $L_p(0,1),\ 1<p<+\infty$ and its subspaces) and to the orders $\alpha, \ 0< \alpha \le 1$ of the fractional derivatives (the case $n-1 < \alpha \le n \in \N$ can be covered by analogy to the known theory of the Riemann-Liouville fractional derivatives of order $\alpha,\ \alpha \in \R_+$).
\end{Rem}

%%%%%%%%%%%%%%%%%%%%%%%%%%%%%%%%%%% 3.1 %% 

\subsection{The Riemann-Liouville fractional derivative}
\label{s_RL}

In this subsection, some known results (see, e.g., \cite{Samko}) are presented in a slightly different  form suitable for our further constructions.

\medskip

We start with the formula \eqref{ftFC} and rewrite it in equivalent form of two equations:
\vskip -14pt
\begin{equation}
\label{ftFC_1}
(I^\alpha\ \phi)(x) = f(x),\ \ (D^\alpha\, f)(x) = \phi(x),\ \,  x\in [0,\, 1].
\end{equation}
\vskip -3pt \noindent
The second of equations from \eqref{ftFC_1} defines the fractional derivative $D^\alpha$ of a function $f$ as the solution $\phi$ of the Abel integral equation with the right-hand side $f$. Now we recall Theorem \ref{t_A} and reformulate it as follows:

\vspace*{-2pt}

\begin{The}
\label{t_A-1}
On the space of functions $X^0 = I^\alpha(L_1(0,1))$, the unique fractional derivative $D^\alpha$ of order $\alpha$, $0< \alpha < 1$ is given by the formula
\vskip -12pt
\begin{equation}
\label{RLDB}
(D^\alpha\, f)(x)  = \frac{d}{dx}\, (I^{1-\alpha}\, f)(x) = \frac{d}{dx}\, \frac{1}{\Gamma(1-\alpha)} \int_0^x (x-t)^{-\alpha}\, f(t)\, dt.
\end{equation}
\end{The}

To justify Theorem \ref{t_A-1}, we just mention that for any function $f\in X^0$ the conditions \eqref{Acond} of Theorem \ref{t_A} are satisfied. This follows from the representation
\begin{equation}
\label{X0}
(I^{1-\alpha}\, f)(x) = (I^{1-\alpha}\, I^\alpha\, \phi)(x) = (I^1 \phi)(x), \phi \in L_1(0,1)
\end{equation}
and the formula  \eqref{absc}. Then the Abel integral equation (the first formula in \eqref{ftFC_1}) has a unique solution given by the formula \eqref{RLDB}.

The first part of the formula \eqref{RLDB} can be used to define the fractional derivative $D^\alpha $ of the order $\alpha = 1$ as the first order derivative:
\begin{equation}
\label{RLDB_0}
(D^1\, f)(x)  = \frac{d}{dx}\, (I^{0}\, f)(x) = \frac{df}{dx}.
\end{equation}

In what follows, we refer to the operator $D^\alpha = \frac{d}{dx}\, I^{1-\alpha}: X^0\to L_1(0,1)$ with $0< \alpha \le 1$ as to the basic Riemann-Liouville fractional derivative
of order $\alpha$ (the term ``basic" refers to the domain $X^0$ of $D^\alpha$).

\vspace*{-2pt}

\begin{Rem}
\label{r_4}
The basic Riemann-Liouville fractional derivative $D^\alpha: X^0\to L_1(0,1)$ is a one-to-one mapping from $X^0$ onto $L_1(0,1)$. For $0<\alpha <1$, this follows from Theorem \ref{t_A}, representation \eqref{ftFC_1}, and Remark \ref{r_2}, whereas for $\alpha=1$ this statement can be easily directly verified.
\end{Rem}

Now let us continue with defining the Riemann-Liouville fractional derivative on the interval $[0,\, 1]$. Evidently, the formula \eqref{RLDB} makes sense for a space of functions larger than $X^0$, namely, for the space
\begin{equation}
\label{X1}
X^1_{RL} = \{ f:\, I^{1-\alpha}\, f \in \mbox{AC}([0,\, 1])\}.
\end{equation}
Indeed,  for $f\in X^1_{RL}$, the representation
\vskip -10pt
\begin{equation}
\label{absc_1}
(I^{1-\alpha} f)(x) = (I^{1-\alpha} f)(0) + \int_0^x \phi(t)\, dt,\ x\in [0,\, 1],\ \phi \in L_1(0,1)
\end{equation}
\vskip -2pt \noindent
holds true (see the formula \eqref{absc}) and thus
\vskip -10pt
\begin{equation}
\label{absc_2}
(D^\alpha\, f)(x) = \frac{d}{dx}(I^{1-\alpha} f)(x) = \phi(x), \ x\in [0,\, 1]
\end{equation}
\vskip -2pt \noindent
according to Remark \ref{r_0}.

\vspace*{-2pt}

\begin{Def}
\label{d_RL}
The extension of the basic Riemann-Liouville fractional derivative $D^\alpha: X^0\to L_1(0,1)$ to the domain $X^1_{RL}$ is called the Riemann-Liouville fractional derivative of order $\alpha,\ 0< \alpha \le 1$:
\vskip -10pt
\begin{equation}
\label{RLD}
(D^\alpha_{RL}\, f)(x) = \frac{d}{dx}(I^{1-\alpha} f)(x),\ D^\alpha_{RL}:\,X^1_{RL}\to L_1(0,1).
\end{equation}
\end{Def}

In contrast to the basic Riemann-Liouville fractional derivative,  the Riemann-Liouville fractional derivative is not injective and its kernel is a one-dimensional vector space:
\vskip -10pt
\begin{equation}
\label{RLD_K}
\mbox{Ker}(D^\alpha_{RL}) =\left\{c_1\, x^{\alpha -1},\ c_1 \in \R\right\}.
\end{equation}
\vskip -2pt \noindent
This immediately follows from the formula \eqref{power} and Remark \ref{r_2}. It is worth mentioning that the basis function $f_1(x)=x^{\alpha -1}$ does not belong to the space $X^0$  because  $(I^{1-\alpha} t^{\alpha -1})(x) \equiv \Gamma(\alpha)\, \forall x\in [0,\ 1]$ that contradicts to the condition $(I^{1-\alpha} f)(0) = 0$ that is fulfilled for any $f\in X^0$ (see the formula \eqref{X0}). Otherwise, we have the inclusions
\vskip -10pt
\begin{equation}
\label{incl1}
X^0 \subset X^1_{RL},\ \mbox{AC}([0,\, 1])\subset X^1_{RL}.
\end{equation}
\vskip -2pt \noindent
The first inclusion follows from the formula \eqref{X0}
and a proof of the second inclusion can be found in \cite{Samko} (Lemma 2.1).  For other properties of the Riemann-Liouville fractional derivative introduced above we refer the readers to \cite{Samko}.

Here, we just mention that for the Riemann-Liouville fractional derivative, the Fundamental Theorem of FC (Theorem \ref{t2}) is evidently valid on even larger space of functions $X^2_{RL} = L_1(0,1),\ X^1_{RL}\subset X^2_{RL}$, i.e., the formula
\vskip -13pt
\begin{equation}
\label{ftFC_RL}
(D^\alpha_{RL}\, I^\alpha\, f)(x) = f(x),\ x\in [0,\, 1],\ f \in X^2_{RL}
\end{equation}
\vskip -3pt \noindent
holds true.

\subsection{The Caputo fractional derivative} %%%%% 3.2 %%%%%%%%%%%%%%%%%%
\label{s_C}

As stated in Theorem \ref{t_A-1}, the basic Riemann-Liouville fractional derivative is the unique one-parameter fractional derivative on the function space $X^0 = I^\alpha (L_1(0,1))$. Its extension to the larger space $X^1_{RL}$ that leads to the Riemann-Liouville fractional derivative is also unique. Because in mathematics one usually works with the maximal domains and extensions of formulas valid on these domains, from the mathematical viewpoint, the Riemann-Liouville fractional derivatives can be considered as the only ``right" one-pa\-ra\-me\-ter family of the fractional derivatives defined  on a finite interval. Indeed, in the classical mathematical literature (see \cite{Samko} and many hundreds references therein), mainly this derivative and  some of its modifications have been considered on a finite interval.    Where is then the place for the Caputo fractional derivative?

The trick with its definition is that first one needs to contract the basic space $X^0$ and to introduce a space of functions, where the first order derivative commutes with the Riemann-Liouville fractional integral of order $1-\alpha$ ($0 < \alpha \le 1$):
\vskip -11pt
\begin{equation}
\label{X0C}
X^0_{C} = \left\{ f\in X^0:\, \frac{d}{dx}\, I^{1-\alpha}\, f = I^{1-\alpha}\, \frac{df}{dx} \right\}.
\end{equation}
In particular, the space $X^0_{C}$ contains the functions
$ f  \in \mbox{AC}([0,\, 1]) $ that satisfy the condition $f(0) = 0$. Indeed, according to \eqref{absc}, these functions
can be represented in the form
\vskip -10pt
$$
f(x) =  \int_0^x \phi(t)\, dt,\ x\in [0,\, 1],\ \phi \in L_1(0,1).
$$
Then we have the following chain of relations ($0\le \alpha$):
$$
(I^\alpha\, \frac{df}{dx})(x) = (I^\alpha\, \frac{d}{dx}\, I^1\, \phi)(x) = (I^\alpha\, \phi)(x) % =
$$
$$
=  \frac{d}{dx} (I^1\, I^\alpha\, \phi)(x) = \frac{d}{dx} (I^\alpha\, I^1\, \phi)(x) =
\frac{d}{dx} (I^\alpha\, f)(x).
$$

The basic Caputo fractional derivative of order $\alpha$, $0 < \alpha \le 1$, is introduced as follows:
\vskip -12pt
\begin{equation}
\label{CDB}
(D^\alpha_C\, f)(x)  = (I^{1-\alpha}\, \frac{df}{dx})(x),  \ D^\alpha_C: \, X^0_{C}\to L_1(0,1).
\end{equation}

Of course, the basic Caputo fractional derivative is identical with the basic Riemann-Liouville fractional derivative restricted to the domain $X^0_{C}$ and thus it is nothing new. However, we get a new operator by extension of its domain! The operator \eqref{CDB}  is well defined, say,  on the space $X^1_{C} = \mbox{AC}([0,\, 1])$.

\vspace*{-2pt}

\begin{Def}
\label{d_C}
The extension of the basic Caputo fractional derivative $D^\alpha_C: X^0_{C}\to L_1(0,1)$ to the domain $X^1_{C}$ is called the Caputo fractional derivative of order $\alpha$, $0 < \alpha \le 1$:
\vskip -10pt
\begin{equation}
\label{CD}
(D^\alpha_{C}\, f)(x) = (I^{1-\alpha} \frac{d}{dx}f)(x),\ \, D^\alpha_{C}:\,X^1_{C}\to L_1(0,1).
\end{equation}
\end{Def}

Evidently, the kernel of the Caputo fractional derivative coincides with the kernel of the first order derivative:
\begin{equation}
\label{CD_K}
\mbox{Ker}(D^\alpha_{C}) =\left\{c_1, \ c_1 \in \R\right\}.
\end{equation}

For the functions from $X^1_{C}$, there is a simple connection between the Riemann-Liouville and the Caputo fractional derivatives (Lemma 2.2 in \cite{Samko}):
\vskip -13pt
\begin{equation}
\label{C_RL}
(D^\alpha_{C}\, f)(x) = (D^\alpha_{RL}\, f)(x)-\frac{f(0)}{\Gamma(1-\alpha)}\, x^{-\alpha},\ \  x>0,\ f\in X^1_{C}.
\end{equation}

As in the case of the Riemann-Liouville fractional derivative, the Fundamental Theorem of FC (Theorem \ref{t2}) for the Caputo fractional derivative is valid on even larger space of functions
\begin{equation}
\label{X2}
X_{FT} = \left\{ f:\, I^{\alpha}f\in \mbox{AC}([0,\, 1])\ \mbox{and}\ (I^{\alpha}f)(0) = 0\right\},
\end{equation}
i.e., the formula
\vskip -10pt
\begin{equation}
\label{ftFC_C}
(D^\alpha_{C}\, I^\alpha\, f)(x) = f(x),\ \, x\in [0,\, 1],\ f \in X_{FT}
\end{equation}
holds true. Let us prove it. According to \eqref{absc}, for a function $f$ from
$X_{FT}$, the representation
\vskip -12pt
\begin{equation}
\label{ftc_1}
(I^\alpha\, f)(x) = (I^1\, \phi)(x), \ x\in [0,\, 1]
\end{equation}
holds true with a function $\phi \in L_1(0,1)$. We then get the following chain of equalities:
\vskip-12pt
$$
(D^\alpha_{C}\, I^\alpha\, f)(x) = (I^{1-\alpha}\, \frac{d}{dx}\, I^\alpha\, f)(x) =
(I^{1-\alpha}\, \frac{d}{dx}\, I^1\, \phi)(x)
= (I^{1-\alpha}\, \phi)(x).
$$
Because for $\phi \in L_1(0,1)$, the fractional integral $I^{1-\alpha}\, \phi$ also belongs to $L_1(0,1)$ (Theorem 2.1 in \cite{Diet10}), we can apply the operator $I^\alpha$ to the last formula:
\vskip -12pt
$$
(I^\alpha\, (D^\alpha_{C}\, I^\alpha\, f)) (x) = (I^\alpha\, I^{1-\alpha}\, \phi)(x) = (I^1\, \phi)(x) = (I^\alpha\, f)(x).
$$
The Riemann-Liouville fractional integral is injective (Remark \ref{r_2}) and thus the formula \eqref{ftFC_C} follows from the last relation.

It is worth mentioning that the space $X_{FT}$ defined by \eqref{X2} can be also characterized as follows (Theorem 2.3 in \cite{Samko}):
\begin{equation}
\label{X2_1}
X_{FT}= I^{1-\alpha}(L_1(0,1))\, (\forall f \in X_{FT}\, \exists \phi \in L_1(0,1):\, f(x) = (I^{1-\alpha}\, \phi)(x)).
\end{equation}

For other properties of the Caputo fractional derivative we refer the readers to \cite{Diet10}.

\subsection{The Hilfer fractional derivative} %%%%%%%% 3.3 %%%%%%%%%%%%%%%
\label{s_H}

The third known family of the fractional derivatives defined on a finite interval that fulfills the Fundamental Theorem of FC is the family of the generalized Riemann-Liouville fractional derivatives. They were introduced by Hilfer in \cite{Hil00} and are nowadays referred to as the Hilfer fractional derivatives.

The schema for construction of the Hilfer fractional derivative of order $\alpha,\ 0< \alpha \le 1$ on a finite interval is the same as the one employed for the Caputo fractional derivative. We start by  defining a suitable basic space of functions, where the first order derivative commutes with a certain Riemann-Liouville fractional integral. Let a parameter $\gamma_1 \in \R$ satisfy the conditions
\vskip -13pt
\begin{equation}
\label{gamma1}
0\le \gamma_1 \le 1-\alpha.
\end{equation}
The space of functions for the basic Hilfer fractional derivative is defined as follows:
\vskip -10pt
\begin{equation}
\label{X0H}
X^0_{H} = \left\{ f\in X^0:\, \frac{d}{dx}\, I^{\gamma_1}\, f = I^{\gamma_1}\, \frac{df}{dx} \right\}.
\end{equation}
As in the case of the space $X^0_{C}$ for the Caputo fractional derivative, the space $X^0_{H}$ contains in particular the functions
$ f  \in \mbox{AC}([0,\, 1]) $ that satisfy the condition $f(0) = 0$.

The basic Hilfer fractional derivative of order $\alpha,\ 0<\alpha \le 1$ and type $\gamma_1$, $0\le \gamma_1 \le 1-\alpha$ is introduced as follows:
\begin{equation}
\label{HDB}
(D^{\alpha,\gamma_1}_H\, f)(x)  = (I^{\gamma_1}\, \frac{d}{dx}\, I^{1-\alpha-\gamma_1}\, f)(x),  \ D^{\alpha,\gamma_1}_H: \, X^0_{H}\to L_1(0,1).
\end{equation}

On the space $X^0_{H}$, the basic Hilfer fractional derivative is identical to the basic Riemann-Liouville fractional derivative restricted to the domain $X^0_{H}$:
\vskip -13pt
$$
(D^{\alpha,\gamma_1}_H\, f)(x)  = (I^{\gamma_1}\, \frac{d}{dx}\, I^{1-\alpha-\gamma_1}\, f)(x) = \frac{d}{dx} \,(I^{\gamma_1}\, I^{1-\alpha-\gamma_1}\, f)(x) % =
$$
$$
=  \frac{d}{dx} \,(I^{1-\alpha}\, f)(x) = (D^{\alpha}_{RL}\, f)(x), \ \, f \in X^0_{H}.
$$
However, the domain of the basic Hilfer fractional derivative can be extended to the larger space of functions:
\begin{equation}
\label{X1H}
X^1_{H} = \left\{ f:\, I^{1 - \alpha -\gamma_1}\, f \in \mbox{AC}([0,\, 1])\right\}.
\end{equation}

\vspace*{-2pt}

\begin{Def}
\label{d_H}
The extension of the basic Hilfer fractional derivative $D^\alpha_H: X^0_{H}\to L_1(0,1)$ to the domain $X^1_{H}$ is called the Hilfer fractional derivative of order  $\alpha,\ 0<\alpha \le 1$ and type $\gamma_1$, $0\le \gamma_1 \le 1-\alpha$:
\vskip -10pt
\begin{equation}
\label{HD}
(D^{\alpha,\gamma_1}_H\, f)(x)  = (I^{\gamma_1}\, \frac{d}{dx}\, I^{1-\alpha-\gamma_1}\, f)(x),\ D^{\alpha,\gamma_1}_H:\, X^1_{H}\to L_1(0,1).
\end{equation}
\end{Def}

\begin{Rem}
\label{rH1}
In \cite{Hil00} and subsequent publications \cite{Hil19, Hil, KocLuc19}, another parametrization of the Hilfer fractional derivative has been employed:
\vskip -10pt
$$
(D^{\alpha,\beta}_H\, f)(x)  = (I^{\beta(1-\alpha)}\, \frac{d}{dx}\, I^{(1-\alpha)(1-\beta)}\, f)(x),\ \,
0 < \alpha \le 1, 0 \le \beta \le 1.
$$
Setting $\gamma_1 = \beta(1-\alpha)$ in this formula, we get the formula \eqref{HD}. In this paper, we prefer using the parametrization  \eqref{HD} because it admits a straightforward generalization as we will see in the subsequent subsections.
\end{Rem}

The Fundamental Theorem of FC (Theorem \ref{t2}) for the Hilfer fractional derivative is valid on the  space $X_{FT}$ defined by \eqref{X2}:
\begin{equation}
\label{ftFC_H}
(D^{\alpha,\gamma_1}_H\, I^\alpha\, f)(x) = f(x),\ x\in [0,\, 1],\ f \in X_{FT},\ 0 < \alpha \le 1,\ 0\le \gamma_1 \le 1-\alpha.
\end{equation}

Its proof follows the steps of the proof of the formula \eqref{ftFC_C} for the Caputo fractional derivative. We start with the representation \eqref{ftc_1} and substitute it into the left-hand side of the formula
\eqref{ftFC_H}:
$$
(D^{\alpha,\gamma_1}_H\, I^\alpha\, f)(x) = (I^{\gamma_1}\, \frac{d}{dx}\,  I^{1-\alpha-\gamma_1}\, I^\alpha\, f)(x) =
(I^{\gamma_1}\, \frac{d}{dx}\, I^{1-\alpha-\gamma_1}\, I^1\, \phi)(x) % =
$$
$$
= (I^{\gamma_1}\, \frac{d}{dx}\, I^1\, I^{1-\alpha-\gamma_1}\,  \phi)(x) = (I^{\gamma_1}\, I^{1-\alpha-\gamma_1}\,  \phi)(x)  = (I^{1-\alpha}\, \phi)(x).
$$
The rest of the proof is exactly the same as the proof of the formula \eqref{ftFC_C} for the Caputo fractional derivative we presented in the previous subsection.

\vspace*{-2pt}

\begin{Rem}
\label{rH2}
For each type $\gamma_1$, $0 \le \gamma_1 \le 1- \alpha$, the Hilfer derivatives $D^{\alpha,\gamma_1}_H$ of orders $\alpha,\ 0< \alpha \le 1$ form the one-parameter families of the fractional derivatives in the sense of Definition \ref{def1} (see the formula \eqref{ftFC_H}).  For $\gamma_1 = 0$, this family coincides with the Riemann-Liouville fractional derivatives, while for $\gamma_1 = 1-\alpha$ we get the Caputo fractional derivatives.
\end{Rem}

The kernel of the Hilfer fractional derivative can be easily calculated by employing the formula \eqref{power} (\cite{Hil}):
\begin{equation}
\label{HD_K}
\mbox{Ker}(D^{\alpha,\gamma_1}_H) =\left\{c_1\, x^{\alpha +\gamma_1 -1},\ c_1 \in \R \right\}.
\end{equation}
In agreement with Remark \eqref{rH1},  for $\gamma_1 =0$, it coincides with the kernel \eqref{RLD_K} of the Riemann-Liouville fractional derivative and for $\gamma_1 =1-\alpha$ with the kernel \eqref{CD_K} of the Caputo fractional derivative.

\vskip 2pt

For other properties of the Hilfer fractional derivative we refer the readers to \cite{Hil00, Hil}.

\subsection{The 2nd level fractional derivative} %%%%%% 3.4 %%%%%%%%%%%%
\label{s_2L}

In turns out that the procedure from the previous subsection can be employed to introduce some new families of the one-parameter fractional derivatives in the sense of Definition \ref{def1}. In this subsection, a family of the derivatives that we call the 2nd level fractional derivatives of order $\alpha,\ 0<\alpha \le 1$ and type $(\gamma_1,\ \gamma_2)$ is discussed.

First we define a suitable basic space of functions. Let the parameters $\gamma_1,\ \gamma_2 \in \R$ satisfy the conditions
\begin{equation}
\label{gamma12}
0\le \gamma_1,\ 0\le \gamma_2,\ \alpha+\gamma_1 \le 1,\ \alpha +\gamma_1+\gamma_2 \le 2.
\end{equation}
In this subsection, we always suppose that these conditions are fulfilled.

The space of functions for the basic 2nd level fractional derivative is defined as follows:
\begin{equation}
\label{X02L}
X^0_{2L} = \left\{ f\in X^0:\, \frac{d}{dx}\, I^{\gamma_1}\, f = I^{\gamma_1}\, \frac{df}{dx}\ \mbox{and} \ \frac{d}{dx}\, I^{\gamma_1+\gamma_2}\, f = I^{\gamma_1+\gamma_2}\, \frac{df}{dx} \right\}.
\end{equation}
The space $X^0_{2L}$ contains in particular the functions
$ f  \in \mbox{AC}([0,\, 1]) $ that satisfy the condition $f(0) = 0$. As already mentioned, for such functions, the first order derivative commutes with the Riemann-Liouville integral of any order $\alpha$, $0\le \alpha$.

The basic 2nd level fractional derivative of order $\alpha,\ 0<\alpha \le 1$ and type $(\gamma_1,\, \gamma_2)$  is first introduced on the space of functions $X^0_{2L}$ as follows:
\vskip -10pt
\begin{equation}
\label{2LDB}
(D^{\alpha,(\gamma_1,\gamma_2)}_{2L}\, f)(x)  = (I^{\gamma_1}\, \frac{d}{dx}\, I^{\gamma_2}\, \frac{d}{dx}\, I^{2-\alpha-\gamma_1-\gamma_2}\, f)(x).
\end{equation}

It maps the space $X^0_{2L}$  into  $L_1(0,1)$ and coincides with the basic Riemann-Liouville fractional derivative restricted to the domain $X^0_{2L}$:
$$
(D^{\alpha,(\gamma_1,\gamma_2)}_{2L}\, f)(x)  = (I^{\gamma_1}\, \frac{d}{dx}\, I^{\gamma_2}\, \frac{d}{dx}\, I^{2-\alpha-\gamma_1-\gamma_2}\, f)(x) % =
$$
$$
=  \frac{d}{dx}\, (I^{\gamma_1+\gamma_2}\, \frac{d}{dx}\, I^{2-\alpha-\gamma_1-\gamma_2}\, f)(x) =
\frac{d}{dx}\, \frac{d}{dx}\, (I^{\gamma_1+\gamma_2}\,  I^{2-\alpha-\gamma_1-\gamma_2}\, f)(x) % =
$$
$$
= \frac{d}{dx}\,\frac{d}{dx} \,(I^1\, I^{1-\alpha}\, f)(x) = \frac{d}{dx} \,(I^{1-\alpha}\, f)(x) = (D^{\alpha}_{RL}\, f)(x), \ \,  f \in X^0_{2L}.
$$

However, the domain of the basic 2nd level fractional derivative can be extended to the larger space of functions:
\begin{equation}
\label{X12L}
X^1_{2L} = \{ f: I^{2-\alpha-\gamma_1-\gamma_2} f, \,
I^{\gamma_2} \frac{d}{dx} I^{2-\alpha-\gamma_1-\gamma_2} f\in \mbox{AC}([0, 1]) \}.
\end{equation}

\begin{Def}
\label{d_2L}
The extension of the basic 2nd level  fractional derivative $D^{\alpha,(\gamma_1,\gamma_2)}_{2L}:\, X^0_{2L}\to L_1(0,\, 1)$ to the domain $X^1_{2L}$ is called the 2nd level  fractional derivative of order  $\alpha,\ 0<\alpha \le 1$ and type $(\gamma_1,\ \gamma_2)$:
\vskip -10pt
\begin{equation}
\label{2LD}
(D^{\alpha,(\gamma_1,\gamma_2)}_{2L}\, f)(x)  = (I^{\gamma_1}\, \frac{d}{dx}\, I^{\gamma_2}\, \frac{d}{dx}\, I^{2-\alpha-\gamma_1-\gamma_2}\, f)(x),
\end{equation}
\vspace*{-4pt}
$$
D^{\alpha,(\gamma_1,\gamma_2)}_{2L}:\, X^1_{2L}\to L_1(0,1).
$$
\end{Def}

\begin{Rem}
\label{rL21}
The operator \eqref{2LD} is called the 2nd level fractional derivative because its formula contains two pairs of compositions of the first order derivatives and the Riemann-Liouville fractional integrals. In this sense, the Riemann-Liouville, the Caputo, and the Hilfer fractional derivatives are the 1st level derivatives.
\end{Rem}

To justify the notation ``fractional derivative", we formulate and prove the Fundamental Theorem of FC (Theorem \ref{t2}) for the 2nd level   fractional derivative.

\vspace*{-3pt}

\begin{The}
\label{t_ft_L2}
On the space $X_{FT}$, the 2nd level  fractional derivative is a left-inverse operator to the Riemann-Lioville fractional integral , i.e., the relation
\vskip -10pt
\begin{equation}
\label{ftFC_2L}
(D^{\alpha,(\gamma_1,\gamma_2)}_{2L}\, I^\alpha\, f)(x) = f(x),\ x\in [0,\, 1]
\end{equation}
holds true for any $f$ from the space $X_{FT}$ defined by \eqref{X2}.
\end{The}

The proof of this theorem follows the lines of the one of the Fundamental Theorem  for the Hilfer fractional derivative (formula \eqref{ftFC_H}). First we substitute  the representation \eqref{ftc_1} of a function from $X_{FT}$ into the left-hand side of the formula \eqref{ftFC_2L} and get the following chain of relations:
$$
(D^{\alpha,(\gamma_1,\gamma_2)}_{2L}\, I^\alpha\, f)(x) = (I^{\gamma_1}\, \frac{d}{dx}\, I^{\gamma_2}\, \frac{d}{dx}\, I^{2-\alpha-\gamma_1-\gamma_2}\, I^\alpha\, f)(x) % =
$$
\vspace*{-3pt}
$$
= (I^{\gamma_1}\, \frac{d}{dx}\, I^{\gamma_2}\, \frac{d}{dx}\, I^{2-\alpha-\gamma_1-\gamma_2}\, I^1\, \phi)(x) =
(I^{\gamma_1}\, \frac{d}{dx}\, I^{\gamma_2}\, I^{2-\alpha-\gamma_1-\gamma_2}\,  \phi)(x) % =
$$
\vspace*{-3pt}
$$
= (I^{\gamma_1}\, \frac{d}{dx}\, I^{1}\, I^{1-\alpha-\gamma_1}\,  \phi)(x) = (I^{\gamma_1}\, I^{1-\alpha-\gamma_1}\,  \phi)(x) = (I^{1-\alpha}\, \phi)(x).
$$
The rest of the proof is exactly the same as the proof of the formula \eqref{ftFC_H} for the Hilfer fractional derivative that we presented in the previous subsection.

Thus, for each type $(\gamma_1,\, \gamma_2)$ that satisfies the conditions \eqref{gamma12}, the 2nd level factional derivatives $D^{\alpha,(\gamma_1,\gamma_2)}_{2L}$ of orders $\alpha,\ 0< \alpha \le 1$ form the one-parameter families of the fractional derivatives in the sense of Definition \ref{def1}.

\vspace*{-4pt}

\begin{Rem}
\label{r-h-2L}
For $1 \le \gamma_2$ or $\alpha +\gamma_1+\gamma_2 \le 1 $,   the 2nd level fractional derivative is reduced to the 1st level Hilfer fractional derivative, respectively:
$$
D^{\alpha,(\gamma_1,\gamma_2)}_{2L} = I^{\gamma_1}\, \frac{d}{dx}\, I^1\,I^{\gamma_2-1} \, \frac{d}{dx}\,I^{2-\alpha-\gamma_1-\gamma_2} %  =
$$
$$
= I^{\gamma_1+\gamma_2 -1}\, \frac{d}{dx}\, I^{1-\alpha-(\gamma_1+\gamma_2-1)} = D^{\alpha,\gamma_1+\gamma_2-1}_{H},
$$
$$
D^{\alpha,(\gamma_1,\gamma_2)}_{2L} = I^{\gamma_1}\, \frac{d}{dx}\, I^{\gamma_2} \, \frac{d}{dx}\,I^{1}\, I^{1-\alpha-\gamma_1-\gamma_2}
=  I^{\gamma_1}\, \frac{d}{dx}\, I^{\gamma_2} I^{1-\alpha-\gamma_1-\gamma_2} % =
$$
$$
 = I^{\gamma_1}\, \frac{d}{dx}\,  I^{1-\alpha-\gamma_1} = D^{\alpha,\gamma_1}_{H}.
 $$
\end{Rem}

In the rest of this subsection we determine the kernel of the 2nd level fractional derivative of order $\alpha,\ 0<\alpha \le 1$ and type $(\gamma_1,\, \gamma_2)$ by applying the formula \eqref{power}.
In accordance with Remark \ref{r-h-2L}, we restrict ourselves to the case of the truly 2nd level fractional derivative, i.e., to the case when the conditions
\vskip -13pt
\begin{equation}
\label{cond_add}
\gamma_2 <1,\ 1<\alpha+\gamma_1+\gamma_2
\end{equation}
hold true.
These  conditions ensure that the kernel is two-dimensional:
\begin{equation}
\label{2LD_K}
\mbox{Ker}(D^{\alpha,(\gamma_1,\gamma_2)}_{2L}) =\left\{c_1 x^{\alpha +\gamma_1 -1}+c_2 x^{\alpha +\gamma_1 +\gamma_2-2},\ c_1,\, c_2 \in \R \right\}.
\end{equation}
The exponents $\sigma_1 = \alpha +\gamma_1 -1$ and $\sigma_2 = \alpha +\gamma_1 +\gamma_2-2$ of the basis functions of the kernel fulfill the inequalities $-1 <\sigma_k\le 0,\ k=1,2$  because of the conditions \eqref{gamma12} and \eqref{cond_add}.

As an example, let us consider the case of  the 2nd level fractional derivative of order $\alpha$, $0< \alpha < 1$ with $\gamma_1 = \gamma_2 = 1-\alpha$. The conditions \eqref{gamma12} and \eqref{cond_add} are evidently satisfied and the fractional derivative takes the form
$$
D^{\alpha,(1-\alpha,1-\alpha)}_{2L}  = I^{1-\alpha}\, \frac{d}{dx}\, I^{1-\alpha}\, \frac{d}{dx}\, I^{\alpha}, \ \,
0< \alpha < 1.
$$
\vskip -3pt \noindent
Its kernel
\vskip -12pt
$$
\mbox{Ker}(D^{\alpha,(1-\alpha,1-\alpha)}_{2L}) =\left\{c_1 +c_2 x^{-\alpha},\ c_1,\, c_2 \in \R \right\}
$$
can be interpreted as a direct sum of the kernels of the Caputo fractional derivative of order $\alpha$ and the Riemann-Liouville fractional derivative of order $1-\alpha$.

\subsection{The $n$th level fractional derivative} %%%%% 3.5 %%%%%%%%
\label{s_nL}

In this subsection, the construction employed in the previous subsection is extended to the case of $n$ compositions of the first order derivatives and appropriate Riemann-Liouville fractional integrals.  As a result, we arrive at infinitely many different families of the one-parameter fractional derivatives that we call the $n$th level fractional derivatives of order $\alpha,\ 0<\alpha \le 1$ and type $\gamma = (\gamma_1,\ \gamma_2,\dots,\gamma_n)$.

In what follows, we suppose that  the parameters $\gamma_1,\ \gamma_2,\dots,\gamma_n \in \R$ satisfy the following conditions:
\vskip-10pt
\begin{equation}
\label{gamma}
0\le \gamma_k \ \mbox{and}\  \alpha + s_k \le k, \ \  k=1,2,\dots,n,
\end{equation}
where, for convenience, we used the notation
\vskip -10pt
\begin{equation}
\label{not}
s_k:= \sum_{i=1}^k\, \gamma_i,\ k= 1,2,\dots,n.
\end{equation}

A suitable space of functions for the basic $n$th level fractional derivative is defined as follows:
\begin{equation}
\label{X0nL}
X^0_{nL} = \left\{ f\in X^0:\, \frac{d}{dx}\, I^{s_k}\, f = I^{s_k}\, \frac{df}{dx},\ \,  k=1,2,\dots,n \right\}.
\end{equation}
The space $X^0_{nL}$ contains in particular the functions
$ f  \in \mbox{AC}([0,\, 1]) $ that satisfy the condition $f(0) = 0$ because for these functions, the first order derivative commutes with the Riemann-Liouville integral of any order $\alpha$, $0\le \alpha$.

On the space of functions $X^0_{nL}$, the basic $n$th level fractional derivative of order $\alpha,\ 0<\alpha \le 1$ and type $\gamma=(\gamma_1,\, \gamma_2,\dots,\gamma_n)$  is introduced as follows:
\begin{equation}
\label{nLDB}
(D^{\alpha,(\gamma)}_{nL}\, f)(x)  = \left(\prod_{k=1}^n (I^{\gamma_k}\, \frac{d}{dx})\right)\, (I^{n-\alpha-s_n}\, f)(x).
\end{equation}

For $n=2$, \eqref{nLDB} is reduced to the 2nd level fractional derivative \eqref{2LDB}.  The operator $D^{\alpha,(\gamma)}_{nL}$ maps the space $X^0_{nL}$  into  $L_1(0,1)$ and is identical with the basic Riemann-Liouville fractional derivative restricted to the domain $X^0_{nL}$:
\vskip -13pt
$$
(D^{\alpha,(\gamma)}_{nL}\, f)(x)  = \left(\prod_{k=1}^n (I^{\gamma_k}\, \frac{d}{dx})\right)\, (I^{n-\alpha-s_n}\, f)(x) % =
$$
$$
= \frac{d}{dx}\, I^{\gamma_1}\, \left( \prod_{k=2}^n (I^{\gamma_k}\, \frac{d}{dx})\right)\, (I^{n-\alpha-s_n}\, f)(x) % =
$$
$$
= \frac{d}{dx}\, \frac{d}{dx}\, I^{\gamma_1+\gamma_2}\, \left(\prod_{k=3}^n (I^{\gamma_k}\, \frac{d}{dx})\right)\, (I^{n-\alpha-s_n}\, f)(x) % =
$$
$$
= \dots =  \left(\frac{d}{dx}\right)^n \,(I^{s_n}\, I^{n-\alpha-s_n}\, f)(x) = \left(\frac{d}{dx}\right)^n \, (I^{n-\alpha}\, f)(x) % =
$$
$$
= \frac{d}{dx} \, (I^{1-\alpha}\, f)(x) =(D^{\alpha}_{RL}\, f)(x), \ \  f \in X^0_{nL}.
$$

As in the case of the 2nd level fractional derivative, the domain of the basic $n$th level  fractional derivative can be extended to a larger space of functions (an empty product is interpreted as the identity operator):
\vskip -10pt
\begin{equation}
\label{X1nL}
X^1_{nL} = \{ f: \left(\prod_{k=i}^n (I^{\gamma_k}\, \frac{d}{dx})\right)\, I^{n-\alpha-s_n}\, f \in \mbox{AC}([0,\, 1]),\ \ i=2,\dots n+1 \}.
\end{equation}

\vspace*{-2pt}

\begin{Def}
\label{d_nL}
The extension of the basic $n$th level  fractional derivative $D^{\alpha,(\gamma)}_{nL}:\, X^0_{nL}\to L_1(0,\, 1)$ to the domain $X^1_{nL}$ is called the $n$th level  fractional derivative of order  $\alpha,\ 0<\alpha \le 1$ and type $\gamma = (\gamma_1,\ \gamma_2,\dots,\gamma_n)$:
\begin{equation}
\label{nLD}
(D^{\alpha,(\gamma)}_{nL}\, f)(x)  = \left(\prod_{k=1}^n (I^{\gamma_k}\, \frac{d}{dx})\right)\, (I^{n-\alpha-s_n}\, f)(x),\ D^{\alpha,(\gamma)}_{nL}:\, X^1_{nL}\to L_1(0,1).
\end{equation}
\end{Def}

Now we proceed with the Fundamental Theorem of FC (Theorem \ref{t2}) for the $n$th level  fractional derivative.

\vspace*{-3pt}

\begin{The}
\label{t_ft_Ln}
On the space $X_{FT}$ defined by \eqref{X2}, the $n$th level  fractional derivative is a left-inverse operator to the Riemann-Lioville fractional integral, i.e., the relation
\begin{equation}
\label{ftFC_nL}
(D^{\alpha,(\gamma)}_{nL}\, I^\alpha\, f)(x) = f(x),\ \  x\in [0,\, 1]
\end{equation}
holds true for any $f$ from the space $X_{FT}$.
\end{The}

The proof of this theorem repeats the arguments presented in the previous subsection for the 2nd level fractional derivative. Here we restrict ourselves to the only place in the proof that looks slightly different compared to the proof of the formula \eqref{ftFC_2L}. Substitution  of the representation \eqref{ftc_1}  into the left-hand side of the formula \eqref{ftFC_nL} leads to the following chain of relations:
\vskip -10pt
$$
(D^{\alpha,(\gamma)}_{nL}\, I^\alpha\, f)(x) =  \left(\left(\prod_{k=1}^n (I^{\gamma_k}\, \frac{d}{dx})\right)\, I^{n-\alpha-s_n}\, I^\alpha\, f\right)(x) % =
$$
\vspace*{-3pt}
$$
=\, \left(\,\left(\prod_{k=1}^n (I^{\gamma_k}\, \frac{d}{dx})\right)\, I^{n-\alpha-s_n}\, I^{1}\, \phi\right)(x)
\,=\,
\left(\,\left(\prod_{k=1}^{n-1} (I^{\gamma_k}\, \frac{d}{dx})\right)\, I^{\gamma_n}\, I^{n-\alpha-s_{n}}\, \phi\right)(x) $$
\vspace*{-3pt}
$$
= \left(\left(\prod_{k=1}^{n-1} (I^{\gamma_k}\, \frac{d}{dx})\right)\, I^{n-1-\alpha-s_{n-1}}\, I^{1}\, \phi\right)(x) = \dots
 = (I^{\gamma_1}\, \frac{d}{dx}\, I^{1-\alpha-s_1}\, I^{1}\, \phi)(x) 
 $$
 \vspace*{-3pt}
 $$
 = (I^{\gamma_1}\, I^{1-\alpha-\gamma_1}\,  \phi)(x) = (I^{1-\alpha}\, \phi)(x).
$$
The rest of the proof is exactly the same as the proof of the formula \eqref{ftFC_2L} for the 2nd level fractional derivative that we presented in the previous subsection.

Thus, for each type $\gamma=(\gamma_1,\, \gamma_2,\dots,\gamma_n)$ satisfying the conditions \eqref{gamma}, the $n$th level fractional derivatives $D^{\alpha,(\gamma)}_{nL}$ of orders $\alpha,\ 0< \alpha \le 1$ form the one-parameter families of the fractional derivatives in the sense of Definition \ref{def1}.

\vspace*{-3pt}

\begin{Rem}
In \cite{DN}, uniqueness and existence of solutions to some Cauchy problems for the fractional differential equations with the operators similar to the $n$th level fractional derivatives $D^{\alpha,(\gamma)}_{nL}$ (in other notations and with other restrictions on the parameters) were considered. However, no connection to the Riemann-Liouville fractional integrals in form of the Fundamental Theorem \ref{t_ft_Ln} was discussed there.
\end{Rem}

\vspace*{-15pt}

\begin{Rem}
\label{rnn-1}
As in the case of the 2nd level fractional derivative (see Remark \eqref{rL21}), the $n$th level fractional derivatives are reduced to the fractional derivatives of the level less than $n$ if  some of the parameters $\gamma_k,\ k=2,\dots,n$ are equal to or grater than one or if the inequality $\alpha+s_n\le n-1$ holds valid.
\end{Rem}

\vspace*{-3pt}

In the rest of this subsection, we consider the truly $n$th level fractional derivatives and suppose that the conditions
\vskip -10pt
\begin{equation}
\label{cond_add_n}
n-1 <\alpha+s_n \ \ \mbox{and}\  \ \gamma_k <1,\ \  k= 2,\dots,n
\end{equation}
\vskip -3pt \noindent
are satisfied (see Remark \ref{rnn-1}).

To determine the kernel of the $n$th level fractional derivative, we again apply the formula \eqref{power}.  Under the conditions \eqref{cond_add_n}, the kernel is $n$-dimensional:
\vskip -13pt
\begin{equation}
\label{nLD_K}
\mbox{Ker}(D^{\alpha,(\gamma)}_{nL}) =\left\{\sum_{k=1}^n c_k x^{\alpha +s_k -k},\ c_k\in \R \right\}.
\end{equation}
The exponents $\sigma_k = \alpha +s_k - k$ of the basis functions of the kernel fulfill the inequalities $-1 <\sigma_k\le 0,\ k=1,2,\dots,n$  because of the conditions \eqref{gamma} and \eqref{cond_add_n}.

In the case, one or several of the conditions from \eqref{cond_add_n} do not hold true, the $n$th level fractional derivatives degenerate to the derivatives of the level less than $n$ and thus their kernels have dimensions less than $n$.

As an example, we consider the case of  the truly $n$th level fractional derivative of order $\alpha$, $0< \alpha \le 1/(n-1)$ and type $(\gamma_1\dots,\gamma_n)$ with $\gamma_k = 1-\alpha, k=1,2,\dots,n$. In this case,  the fractional derivative takes the form
$$
D^{\alpha,(1-\alpha)}_{nL}  = I^{1-\alpha}\, \frac{d}{dx} \dots I^{1-\alpha}\frac{d}{dx}\, I^{(n-1)\alpha}.
$$
For $0< \alpha \le 1/(n-1)$, the conditions \eqref{gamma} and \eqref{cond_add_n} are satisfied and the kernel of this fractional derivative is $n$-dimensional:
$$
\mbox{Ker}(D^{\alpha,(1-\alpha)}_{2L}) =\left\{\sum_{k=1}^n c_k x^{-\alpha(k-1)},\ c_k\in \R  \right\}.
$$
It can be interpreted as a direct sum of the kernels of the Caputo fractional derivative of order $\alpha$ and the Riemann-Liouville fractional derivatives of orders $1-\alpha\, k$,\, $k=1,\dots,n-1$.

%%%%%%%%%%%%%%%% Section 4 %%%%%%%%%%%%%%%%%%%%%%%%%%%%%%%%%%%

\section{Some properties of the fractional integrals and derivatives}
\label{sec:4}

\setcounter{section}{4}
\setcounter{The}{0}

In this section, we mainly address the 2nd level fractional derivatives defined by \eqref{2LD} and their connection to the Riemann-Liouville fractional integral. On the one hand, this derivative is a new object that generalizes the Hilfer fractional derivative that in its turn contains the Riemann-Liouville and the Caputo fractional derivatives as its particular cases. On the other hand, as we have seen in the previous section, the results obtained for the 2nd level fractional derivative can be easily extended to the case of the $n$th level fractional derivatives.

\subsection{Projector of the 2nd level fractional derivative}

The projector  $P^\alpha_{2L}$ of the 2nd level fractional derivative \eqref{2LD} is defined as follows:
\begin{equation}
\label{pro}
(P^\alpha_{2L}\, f)(x) = (\mbox{Id}-I^{\alpha} D^{\alpha,(\gamma_1,\gamma_2)}_{2L}\, f)(x).
\end{equation}

An explicit formula for the projector  is central for one of the most important methods of analysis of the fractional differential equations, namely, for reduction of the fractional differential equations to certain integral equations of Volterra type. The coefficients in the representation formulas for the projectors determine the form of the ``natural" initial conditions required for the corresponding fractional differential equations.

\vspace*{-3pt}

\begin{The}
\label{t-pro}
Under the conditions \eqref{cond_add},  the projector \eqref{pro} of the 2nd level fractional derivative \eqref{2LD} for the functions from the space $X^1_{2L}$ defined by \eqref{X12L} has the following form:
\begin{equation}
\label{pro1}
(P^\alpha_{2L}\, f)(x) = p_1\, x^{\alpha +\gamma_1-1} + p_2\, x^{\alpha +\gamma_1+\gamma_2-2},
\end{equation}
\begin{equation}
\label{p1}
p_1 = \frac{1}{\Gamma(\alpha +\gamma_1)}\left( I^{\gamma_2}\, \frac{d}{dx} \, I^{2-\alpha-\gamma_1-\gamma_2}\, f\right)(0),
\end{equation}
\begin{equation}
\label{p2}
p_2= \frac{1}{\Gamma(\alpha +\gamma_1+\gamma_2 -1)}\left(  I^{2-\alpha-\gamma_1-\gamma_2}\, f\right)(0).
\end{equation}
\end{The}

We prove the theorem by applying a method that works also in the case of the $n$th level fractional derivative. Let us introduce an auxiliary function
\begin{equation}
\label{g}
g(x) := (I^{\alpha} D^{\alpha,(\gamma_1,\gamma_2)}_{2L}\, f)(x).
\end{equation}
For $f\in X^1_{2L}$, the derivative $D^{\alpha,(\gamma_1,\gamma_2)}_{2L}\, f$ is well defined and belongs to $L_1(0,1)$ (see subsection \ref{s_2L}) and thus the function $g$ is from the space $I^\alpha(L_1(0,1))$. Then we can act on $g$ with  the operator $D^{\alpha,(\gamma_1,\gamma_2)}_{2L}$ and apply the Fundamental Theorem of FC for the 2nd level fractional derivative:
$$
(D^{\alpha,(\gamma_1,\gamma_2)}_{2L}\, g)(x) = (D^{\alpha,(\gamma_1,\gamma_2)}_{2L}\, I^{\alpha} D^{\alpha,(\gamma_1,\gamma_2)}_{2L}\, f)(x) = (D^{\alpha,(\gamma_1,\gamma_2)}_{2L}\, f)(x).
$$
Thus the function $g -f$ belongs to the kernel of $D^{\alpha,(\gamma_1,\gamma_2)}_{2L}$ given by the formula \eqref{2LD_K} and we arrive at the representation
\begin{equation}
\label{f-g}
g(x) = f(x) + c_1 x^{\alpha +\gamma_1 -1}+c_2 x^{\alpha +\gamma_1 +\gamma_2-2}.
\end{equation}
To determine the coefficients $c_1,\, c_2$, we first apply the operator $I^{1 -\alpha-\gamma_1}$ to the function $g$:
\vskip -10pt
$$
(I^{1-\alpha-\gamma_1 }\, g)(x) = (I^{1-\alpha-\gamma_1 } I^{\alpha} D^{\alpha,(\gamma_1,\gamma_2)}_{2L}\, f)(x) % =
$$
$$
= \left(I^{1-\alpha-\gamma_1 } I^{\alpha} I^{\gamma_1}\, \frac{d}{dx}\, I^{\gamma_2}\, \frac{d}{dx}\, I^{2-\alpha-\gamma_1-\gamma_2}\, f\right)(x) % =
$$
$$
= \left(I^{1}\, \frac{d}{dx}\, I^{\gamma_2}\, \frac{d}{dx}\, I^{2-\alpha-\gamma_1-\gamma_2}\, f\right)(x) % =
$$
$$
= \left(I^{\gamma_2}\, \frac{d}{dx}\, I^{2-\alpha-\gamma_1-\gamma_2}\, f\right)(x) - \left(I^{\gamma_2}\, \frac{d}{dx}\, I^{2-\alpha-\gamma_1-\gamma_2}\, f\right)(0).
$$
Then we apply the operator  $I^{1-\gamma_2}$ to the result of the previous evaluation:
$$
(I^{1-\gamma_2}\, I^{1 -\alpha-\gamma_1}\, g)(x) = \left(I^{1}\, \frac{d}{dx}\, I^{2-\alpha-\gamma_1-\gamma_2}\, f\right)(x) %% -
$$
$$
- (I^{1-\gamma_2}\, \left(I^{\gamma_2}\, \frac{d}{dx}\, I^{2-\alpha-\gamma_1-\gamma_2}\, f\right)(0))(x) % =
$$
$$
= \left(I^{2-\alpha-\gamma_1-\gamma_2}\, f\right)(x) - \left(I^{2-\alpha-\gamma_1-\gamma_2}\, f\right)(0) %% -
$$
$$
- \frac{1}{\Gamma(2-\gamma_2)}\left(I^{\gamma_2}\, \frac{d}{dx}\, I^{2-\alpha-\gamma_1-\gamma_2}\, f\right)(0)\, x^{1-\gamma_2}.
$$
On the other hand, we can apply the operator $I^{1-\alpha-\gamma_1 }\, I^{1-\gamma_2} = I^{2-\alpha-\gamma_1 -\gamma_2}$ to the representation \eqref{f-g} and get the relation
$$
(I^{2-\alpha-\gamma_1 -\gamma_2}\, g)(x) = (I^{2-\alpha-\gamma_1 -\gamma_2}\, f)(x) %% +
$$
$$
+ \, c_1\frac{\Gamma(\alpha +\gamma_1)}
{\Gamma(2-\gamma_2)} x^{1-\gamma_2}+c_2 \frac{\Gamma(\alpha +\gamma_1 +\gamma_2 -1)}{\Gamma(1)}x^{0}.
$$
Comparing the coefficients by the same powers of $x$ in the last two formulas, we arrive at the values of the coefficients $c_1$ and $c_2$:
$$
c_1 = -\frac{1}{\Gamma(\alpha +\gamma_1)}\left( I^{\gamma_2}\, \frac{d}{dx} \, I^{2-\alpha-\gamma_1-\gamma_2}\, f\right)(0),
$$
\vspace*{-3pt}
$$
c_2= -\frac{1}{\Gamma(\alpha +\gamma_1+\gamma_2 -1)}\left(  I^{2-\alpha-\gamma_1-\gamma_2}\, f\right)(0).
$$
The statement of the theorem follows now from these formulas and the representation \eqref{f-g}.

\begin{Rem}
\label{r_2ndFT}
Theorem \ref{t-pro} can be rewritten in form of the 2nd Fundamental Theorem of FC for the 2nd level fractional derivative:
\begin{equation}
\label{2ndFT}
(I^{\alpha} D^{\alpha,(\gamma_1,\gamma_2)}_{2L}\, f)(x) = f(x) - p_1\, x^{\alpha +\gamma_1-1} - p_2\, x^{\alpha +\gamma_1+\gamma_2-2},
\end{equation}
where the coefficients $p_1$ and $p_2$ are defined as in \eqref{p1} and \eqref{p2}, respectively.
\end{Rem}

\vspace*{-15pt}

\begin{Rem}
As mentioned in the previous section, if one of the conditions \eqref{cond_add} does not hold true, the 2nd level fractional derivatives are reduced to the Hilfer fractional derivatives and their kernels become one-dimensional. In these cases, one of the coefficients $c_1$ or $c_2$ in the representation \eqref{f-g} and thus one of the coefficients $p_1$ or $p_2$ in the formula \eqref{pro1} for the projector
$P^\alpha_{2L}$ is equal to zero.
\end{Rem}

Say, for $\gamma_2 = 1$ or $\alpha +\gamma_1+\gamma_2\le 1$, the 2nd level fractional derivative $D^{\alpha,(\gamma_1,\gamma_2)}_{2L}$ is reduced to the Hilfer fractional derivative $D^{\alpha,\gamma_1}_H$ (see Remark \ref{rL21}). As a result, the projector $P_H^\alpha$ of the Hilfer fractional derivative  takes the following known form (\cite{Hil}):
\vskip -10pt
\begin{equation}
\label{pro1_H_1}
(P^\alpha_{H}\, f)(x) = \frac{1}{\Gamma(\alpha +\gamma_1)}\left(  I^{1-\alpha-\gamma_1}\, f\right)(0)\, x^{\alpha +\gamma_1-1}.
\end{equation}
Substituting $\gamma_1 = 0$ into the formula \eqref{pro1_H_1}, we get the projector of the Riemann-Liouville fractional derivative
\begin{equation}
\label{pro1_RL}
(P^\alpha_{RL}\, f)(x) = \frac{1}{\Gamma(\alpha)}\left(  I^{1-\alpha}\, f\right)(0)\, x^{\alpha -1}.
\end{equation}
The value $\gamma_1 = 1-\alpha$ corresponds to the Caputo fractional derivative:
\begin{equation}
\label{pro1_C}
(P^\alpha_{C}\, f)(x) = f(0).
\end{equation}
The last formula makes clear why many researchers prefer to employ the Caputo fractional derivative while working with the fractional differential equations: it is the only fractional derivative on a finite interval that admits the same initial conditions as the ones usually posed for the differential equations with the integer order derivatives.

\subsection{Laplace transform of the 2nd level fractional derivative} %%%% 4.2 %%%%%%%

The theory of the fractional derivatives on a finite interval that we addressed until now can be extended to the case of the semi-axis following the approach employed for the Riemann-Liouville fractional derivative (\cite{Samko}). This will be done elsewhere. Here, we just introduce the 2nd level fractional derivative on the positive real semi-axis and derive a formula for its Laplace transform.

For the functions from the space $L_{loc}(\R_+)$, the Riemann-Liouville fractional integral on the positive real semi-axis is defined as follows:
\begin{equation}
\label{RLI_R}
(I^\alpha_{0+}\, f)(x) = \begin{cases}
\frac{1}{\Gamma(\alpha)} \int_0^x (x-t)^{\alpha -1}\, f(t)\, dt,\ x>0,& \alpha >0,\\
f(x),\ x>0, & \alpha = 0.
\end{cases}
\end{equation}

For $\alpha >0$, the Riemann-Liouville fractional integral $I^\alpha_{0+}$ can be interpreted as the Laplace convolution of the functions $f=f(x)$ and $h_\alpha(x) = x^{\alpha -1}/\Gamma(\alpha)$, $x>0$. Application of the convolution theorem for the Laplace transform leads to the well-known result  for the Laplace transform of the Riemann-Liouville fractional integral $I^\alpha_{0+}$ (\cite{Samko})
\begin{equation}
\label{RLI_L}
({\mathcal  L}\, I^\alpha_{0+}\, f)(s) = s^{-\alpha}\, ({\mathcal  L}\, f)(s),\ \Re(s) > \max\{s_f,\, 0\}
\end{equation}
that is valid under the condition that the Laplace transform of the function $f$ given by the integral
\vskip -12pt
\begin{equation}
\label{Lap}
({\mathcal  L}\, f)(s) = \int_0^{+\infty} f(t)\, e^{-st}\, dt
\end{equation}
\vskip -3pt \noindent
does exist for $\Re(s) > s_f$.

On the positive real semi-axis, the 2nd level fractional derivative of order $\alpha,\ 0<\alpha \le 1$ and type $(\gamma_1,\, \gamma_2)$ is defined as follows:
\begin{equation}
\label{2LDB_R}
(D^{\alpha,(\gamma_1,\gamma_2)}_{2L_+}\, f)(x)  = (I^{\gamma_1}_{0+}\, \frac{d}{dx}\, I^{\gamma_2}_{0+}\, \frac{d}{dx}\, I^{2-\alpha-\gamma_1-\gamma_2}_{0+}\, f)(x).
\end{equation}

For the appropriate spaces of functions, both the Fundamental Theorem of FC and the projector formula that we derived for the  2nd level fractional derivative on a finite interval are valid for the operators $I^{\alpha}_{0+}$ and $D^{\alpha,(\gamma_1,\gamma_2)}_{2L_+}$ defined on the positive real semi-axis (the proofs repeat the arguments we applied in the case of a finite interval).

In what follows we restrict ourselves to the truly 2nd level fractional derivatives and suppose that the conditions \eqref{cond_add} are satisfied. Then the 2nd Fundamental Theorem of FC (Remark \ref{r_2ndFT}) holds true and we have the representation
\vskip -10pt
$$
(I^{\alpha}_{0+} D^{\alpha,(\gamma_1,\gamma_2)}_{2L_+}\, f)(x) = f(x) - p_1\, x^{\alpha +\gamma_1-1} - p_2\, x^{\alpha +\gamma_1+\gamma_2-2},
$$
where the coefficients $p_1$ and $p_2$ are given by \eqref{p1} and \eqref{p2}, respectively.

Note that under the conditions \eqref{cond_add}, the exponents $\sigma_1 = \alpha +\gamma_1-1$ and $\sigma_2 = \alpha +\gamma_1+\gamma_2-2$ satisfy the inequalities $-1 <\sigma_k\le 0,\ k=1,2$ and thus we can
apply the Laplace transform to the last formula. Using \eqref{RLI_L}, we get the equation
\vskip -12pt
$$
s^{-\alpha}\, ({\mathcal  L}\,  D^{\alpha,(\gamma_1,\gamma_2)}_{2L_+}\, f)(s) = ({\mathcal  L}\, f)(s) - p_1\, \frac{\Gamma(\alpha+\gamma_1)}{ s^{\alpha + \gamma_1}}  - p_2\,  \frac{\Gamma(\alpha+\gamma_1+\gamma_2 -1)}{ s^{\alpha + \gamma_1+\gamma_2-1}}.
$$
Multiplying it with $s^\alpha$ leads to the  formula
\vskip -10pt
$$
({\mathcal  L}\,  D^{\alpha,(\gamma_1,\gamma_2)}_{2L_+}\, f)(s) = s^{\alpha}\, ({\mathcal  L}\, f)(s) - p_1\, \frac{\Gamma(\alpha+\gamma_1)} {s^{\gamma_1}}  - p_2\,  \frac{\Gamma(\alpha+\gamma_1+\gamma_2 -1)}{s^{\gamma_1+\gamma_2-1}}
$$
that can be transformed to the final form
\begin{equation}
\label{Lap_D}
({\mathcal  L}\,  D^{\alpha,(\gamma_1,\gamma_2)}_{2L_+}\, f)(s) = s^{\alpha}\, ({\mathcal  L}\, f)(s) - a_1\, s^{ - \gamma_1}  - a_2\,   s^{- \gamma_1-\gamma_2+1}
\end{equation}
\vskip -2pt \noindent
with
\vskip -13pt
\begin{equation}
\label{a1}
a_1 = \left( I^{\gamma_2}_{0+}\, \frac{d}{dx} \, I^{2-\alpha-\gamma_1-\gamma_2}_{0+}\, f\right)(0), \
a_2= \left(  I^{2-\alpha-\gamma_1-\gamma_2}_{0+}\, f\right)(0).
\end{equation}

If one of the conditions \eqref{cond_add} does not hold true, the kernel of the 2nd level fractional derivative becomes one-dimensional and one of the coefficients $a_1$ or $a_2$ in \eqref{a1} is equal to zero.  In particular, for the Hilfer fractional derivative $D^{\alpha,\gamma_1}_{H_+}$ defined on the real positive semi-axis ($\gamma_2 = 1$  or $\alpha +\gamma_1+\gamma_2\le 1$ in  \eqref{2LDB_R}), we get the formula
\begin{equation}
\label{Lap_D_H}
({\mathcal  L}\,  D^{\alpha,\gamma_1}_{H_+}\, f)(s) = s^{\alpha}\, ({\mathcal  L}\, f)(s) - \left(  I^{1-\alpha-\gamma_1}_{0+}\, f\right)(0)\,   s^{- \gamma_1}.
\end{equation}
\vskip -2pt \noindent
Setting $\gamma_1$ to zero in \eqref{Lap_D_H} leads to the known formula for the Laplace transform of the Riemann-Liouville fractional derivative
\begin{equation}
\label{Lap_D_RL}
({\mathcal  L}\,  D^{\alpha}_{0+}\, f)(s) = s^{\alpha}\, ({\mathcal  L}\, f)(s) - \left(  I^{1-\alpha}_{0+}\, f\right)(0),
\end{equation}
whereas the value $\gamma_1 = 1- \alpha$ in the formula \eqref{Lap_D_H} corresponds to the Laplace transform of the Caputo fractional derivative:
 \vskip -8pt
\begin{equation}
\label{Lap_D_C}
({\mathcal  L}\,  D^{\alpha}_{C_+}\, f)(s) = s^{\alpha}\, ({\mathcal  L}\, f)(s) - f(0)\,   s^{\alpha -1}.
\end{equation}

\vspace*{-2pt}

\begin{Rem}
\label{r_f}
In this paper, we restricted ourselves to analysis of some families of the one-parameter fractional integrals and derivatives on a finite interval and the positive real semi-axis. Thus, the many parameters families of the fractional integrals and derivatives (say, the
Erd\'{e}lyi-Kober operators \cite{Kir94, Yak-Luc} or their Caputo type modifications \cite{LT07, KirLuc13}), the distributed order fractional derivatives \cite{Che02, Che03, Koc08, KocLuc19}, the fractional derivatives with the general kernels \cite{Koch11, Zac08, Zac09}, and the fractional integrals and derivatives defined on the whole space like the Riesz fractional potentials and derivatives \cite{Kwa, Samko} are not included into this theory.
\end{Rem}

\vspace*{-17pt}

\begin{Rem}
\label{r_s}
Another important challenge in FC is studying the fractional order operators on different spaces of functions that are often chosen to satisfy some special requirements of concrete problems. In this paper, we addressed the fractional integrals and derivatives on the space $L_1(0,1)$ and its subspaces. However, the FC operators have been already considered on $L_p(0,1),\ p>1$ (\cite{Samko}), $C_\alpha(0,\infty),\ \alpha >-1$ (\cite{Dim1, Kir94, Luc, LucYak, Yak-Luc}), and the fractional Sobolev spaces (\cite{GLY15}) to mention only some of the relevant spaces. It is worth to  once again stress  that the properties of the FC operators essentially depend on the spaces of functions, where they are defined.
\end{Rem}

\vspace*{-15pt}
%%%%%%%%%%%%%%%%%%%%%%%%%%%%%%%%%%%%%%%%%%%%%%%%%%%%

\end{document}